\documentstyle[amscd]{amsart}
\numberwithin{equation}{section}
\theoremstyle{plain}
\newtheorem{thm}{Theorem}[section]
\newtheorem{cor}[thm]{Corollary}
\newtheorem{lem}[thm]{Lemma}
\newtheorem{prop}[thm]{Proposition}
\begin{document}
\title{Zeta functions and topological entropy of the Markov-Dyck shifts }
\author{Wolfgang Krieger}
\address{Institute for Applied  Mathematics,
University of Heidelberg,
Im Neuenheimer Feld 294,
69120 Heidelberg, Germany}
\author{Kengo Matsumoto}
\address{ 
Department of Mathematical Sciences, 
Yokohama City  University,
22-2 Seto, Kanazawa-ku, Yokohama, 236-0027 Japan}
\maketitle
\begin{abstract}
The Markov-Dyck shifts arise from finite directed graphs. 
An expression for the zeta function of a Markov-Dyck shift is given. 
The derivation of this expression is 
based on a formula in Keller (G. Keller, {\it  Circular codes, loop counting, and zeta-functions}, 
J.\ Combinatorial Theory
{\bf  56}
(1991),
pp.\ 75--83). 
For a class of examples that includes the Fibonacci-Dyck shift
the zeta functions and topological 
entropy ae determined.
\end{abstract}

\def\C{{{\cal C}}}
\def\V{{{\cal V}}}
\def\det{{{\operatorname{det}}}}
\def\trace{{{\operatorname{trace}}}}
\def\card{{{\operatorname{card}}}}




\bigskip

\section{Introduction}
Let $\Sigma$ be a finite alphabet, and let $S_\Sigma $ be the left shift on
$\Sigma^{\Bbb Z}$,
$$
  S_\Sigma((x_i)_{ i \in {\Bbb Z}})   = (x_{i+1})_{ i \in {\Bbb Z}}, \qquad  
(x_i)_{ i \in {\Bbb Z}} \in \Sigma^{\Bbb Z}.
$$
The closed shift-invariant subsystems of the shifts $S_\Sigma $ are called
subshifts. For an introduction to their theory, which belongs to symbolic
dynamics, we refer to \cite{Ki} and \cite{LM}. 
A finite word in the symbols of 
$S_\Sigma$  is called admissible for the subshift 
$X \subset \Sigma^{\Bbb Z} $ if it
appears somewhere in a point of $X$. A subshift is uniquely determined by its
language of admissible words that we denote by ${\cal L}(X)$. 
${\cal L}_n(X)$ 
will denote the set of words in ${\cal L}(X)$ 
of length $n \in {\Bbb N}$.
The topological entropy of the
subshift $X \subset \Sigma^{\Bbb Z}  $
is  given by
$$
h(X) = \lim_{n \to \infty} \frac {1}{n} \log \card {\cal L}_n(X).
$$
Denoting by
$\Pi_n(X)$ the number of points of period $n$ of the subshift $X
\subset \Sigma^{\Bbb Z} ,$
the zeta function of the subshift $X \subset \Sigma^{\Bbb Z}  $ is given by
$$
\zeta_X(z) = e^{\sum_{n \in {\Bbb N}} \frac{\Pi_n(X)z^n}{n}}.
$$
In this paper we are concerned with a class of subshifts that arise from finite
directed graphs as a
special case of constructions that were described in 
\cite{Kr1}, \cite{Kr2}, \cite{HIK}.
Following the line of terminology of \cite{Ma4}, we
call these subshifts Markov-Dyck shifts.
Let $G$ be a finite directed graph with vertex set ${\cal V}$ and edge set 
${\cal E}$.
We denote the initial vertex of $e \in {\cal E}$
by $s(e)$ and the final vertex by $r(e)$.
Let $G^-$ be that graph with vertex set 
${\cal V}$
and 
edge set 
${\cal E}^-$
a copy of ${\cal E}$.
Reverse the directions of the edges in ${\cal E}^-$
to obtain the reversed graph 
$G^+$ of $G^-$
with vertex set ${\cal V}$ and edge set 
${\cal E}^+$.
Denote by ${\cal P}^-(\text{resp. } {\cal P}^+)$
the set of finite paths in $G^- (\text{resp. } G^+)$.
The mapping $e^- \rightarrow e^+ \ ( e^- \in {\cal E}^-)$
extends to the bijection 
$w^- \rightarrow w^+ (w^- \in {\cal P}^-)$  
of 
${\cal P}^-$ onto ${\cal P}^+$
that reverses direction.
With idempotents $P_v, v \in {\cal V},$ 
the set 
${\cal E}^- \cup \{P_v:v \in {\cal V}\}\cup {\cal E}^+$
is the generating set of the graph inverse semigroup of $G$, where, besides
$P_u^2 = P_u, v \in {\cal V},$ the relations are
(see for instance \cite{Pa})
$$
P_uP_w = 0, \qquad   u, w \in {\cal V}, u \ne w,  
$$
\begin{equation}
f^-g^+ 
=
\begin{cases} 
P_{s(f)}, & (f = g),\\
0  &  (s(f) \ne s(g), f, g \in {\cal E} ),\\
\end{cases} 
\end{equation}
and
\begin{align*}
g^+ f^- & = 0, \qquad s(f) \ne s(g),\\
f^- g^-  & = 0, \qquad r(f) \ne s(g), \quad f,g \in {\cal E}. 
\end{align*}
The alphabet of the Markov-Dyck shift $D_G$ of $G$ is
${\cal E}^-\cup {\cal E}^+$
and a word 
$(e_k)_{1 \le k \le K}$ 
is admissible for $D_G$ precisely if 
$$
\prod_{1 \le k \le K}e_k \ne 0.
$$
For the directed graph with one vertex and loops $e_n, 1 \le n \le N, N >1,$
the relations take the form
\begin{equation}
e_n^- e_n^+ =  {\bold 1}, \quad 1 \le n\le N, \qquad  
e_l^- e_m^+ = 0, \quad 1 \le l, m \le N, \ l \ne m.
\end{equation}
and one sees the Dyck inverse monoid \cite{NP}, together with the Dyck shifts 
that were first described in \cite{Kr1}.

The relations (1.2) can be viewed as the multiplicative relations among the
relations that are satisfied by generators of a Leavitt algebra 
\cite{L} or a Cuntz algebra \cite{C} and the relations (1.1) can be viewed as as the multiplicative
relations among the relations that are satisfied by generators of a Leavitt
path algebra of the directed graph G \cite{AA}, or the generators of the
graph $C^*$-algebra of $G$ \cite{CK}, \cite{EW}. 

The zeta
function of the Dyck shifts were determined in \cite{Ke}, and K-theoretic
invariants were computed in \cite{Ma2} and \cite{KM}. 
For related systems, the Motzkin
shifts that add a symbol ${\bold 1}$ to the alphabets of the Dyck shifts,
the zeta functions were determined in \cite{In} and K-theoretic invariants were
computed in \cite{Ma1}.
In section 2 we will obtain an expression for the zeta function of a
Markov-Dyck shift by applying a formula of Keller's \cite{Ke}. 
In section 3 we
derive estimates for the topological entropy of the subsystems of the
Markov-Dyck shifts that are obtained by allowing the paths in the Markov-Dyck
shift to go from ${\cal E}^-$ to ${\cal E}^-$  or vice-versa only when entering a
given vertex. A fortiori, this gives also estimates for the topological
entropy of the Markov Dyck shifts.
In section 4 
we determine the zeta functions and 
topological entropy of the Markov-Dyck shifts 
that arise from directed graphs with  adjacency matrix 
$F(a,b,c) = 
\begin{bmatrix}
a & b \\
c & 0
\end{bmatrix},
\ a,b,c \in {\Bbb N}$
(compare here section 10.3 of \cite{Ma3}).
 K-theoretic
invariants of $D_{F(1,1,1)}$ shift were computed in \cite{Ma5}. 

The length of a word $w$ we denote by $\ell(w)$ and
we denote the generating function of a formal language ${\cal L}$
by 
$g_{\cal L}$,
$$
g_{\cal L}(z)
= \sum_{n \in {\Bbb Z}_+} \card\{ w \in {\cal L} : \ell(w) = n \} z^n.
$$

\section{Zeta functions}


Keller \cite{Ke} has introduced the notion of a circular Markov code.
Here we find ourselves in a situation 
where we will want a Markov code to be given by a set $\C$ of non-empty words
in the symbols of a finite alphabet $\Sigma$
together with a finite set ${\cal V}$ and mappings
$r:\C \rightarrow \V, s:\C \rightarrow \V$.
To
$(\C,r,s)$ there is associated the shift invariant set 
$X_\C \subset \Sigma^{\Bbb Z}$
of points $x \in \Sigma^{\Bbb Z}$ such that 
there are indices 
 $I_k, k \in {\Bbb Z},$ such that 
\begin{equation}
I_0 \le 0 < I_1,\quad  I_k < I_{k+1}, \quad k \in {\Bbb Z}, 
\label{eqn:2.1}
\end{equation}
and such that
\begin{equation}
x_{[I_k,I_{k+1})} \in \C, \qquad k \in {\Bbb Z}, 
\label{eqn:2.2}
\end{equation}
and
\begin{equation}
r(x_{[I_k,I_{k+1})}) = s(x_{[I_{k+1},I_{k+2})}), \qquad k \in {\Bbb Z}. 
\label{eqn:2.3}
\end{equation}
$(\C, r,s)$ is said to be a circular Markov code
if for every periodic point $x$ in $X_\C$ 
the indices $I_k, k \in {\Bbb Z},$ 
such that \eqref{eqn:2.1},  \eqref{eqn:2.2}, 
and  \eqref{eqn:2.3} hold,
are uniquely determined by $x$
and can then be denoted by $I_k(x), k \in {\Bbb Z}$. If ${\cal V}$ contains one
element then one has a circular code (see e.g. \cite{BP}).

Generalizing the formula for the zeta function of $X_\C$, 
where $\C$ is a circular code,
Keller \cite{Ke} has proved a formula for the zeta function 
of $X_\C$, where $\C$ is a circular Markov code.
For completeness we reproduce here Keller's proof 
for the special case that  we we have in mind.

Given a circular Markov code
$(\C, s, r)$ denote by
$
\C(u,w)
$
the set of words $c \in \C$ such that
$s(c) = u$,  $r(c) = w$,$u,w \in {\cal V}.$
Set
\begin{equation*}
g_{\C(u,v), n} = \card \{ c \in \C : s(c) = u, \ r(c) = v, \ \ell(c) =n\}, 
\end{equation*}
and introduce the matrix
$$
H^{(\C)}(z) = (g_{\C(u,v)}(z))_{u,v \in {\cal V}}.
$$

\begin{thm}[Keller]\label{thm:2.1}
For a circular Markov code $(\C, s, r),$ 
$$
\zeta_{X_\C}(z) = \det ( I - H^{(\C)}(z))^{-1}.
$$
\end{thm}
\begin{pf}
Let $n \in {\Bbb N}$.
Consider triples of the form
$(j,c_1, c_1\cdots c_k),$
where
$k \in {\Bbb N},$
and where
\begin{align*}
c_l \in \C, \quad 1 \le l \le k,& \quad
s(c_1) = r(c_k),
r(c_l) = s(c_{l+1}), \quad 1 \le l \le k,\\
&\ell(c_1\cdots c_k) = n - \ell(c_1), \quad j=\ell(c_1)
\end{align*}
To every point $x \in X_\C$ of period $n$ 
one assigns a triple of this kind,
where $k \in {\Bbb N}$ is determined by
$$
n = I_k(x) - I_0(x), \qquad 
j = - I_0(x),
$$
and
$$
c_l = x_{[I_{l-1}(x), I_l(x))}, \quad 1 \le l < k.
$$  
Due to the circularity of 
$(\C,r,s)$
this assignment is bijective.

Denote by $\gamma_i(f)$the $i$-th coefficient of the power series of a
function $f$.
From
\begin{align*}
\gamma_j((H^k)_{u,w})
= & \card \{ c_1\cdots c_k : c_l \in \C, 1 \le l \le k,
s(c_1) = u, r(c_k)= w, \\
 & r(c_l) = s(c_{l+1}), 1 \le l \le k,\ 
\ell(c_1\cdots c_k) = j  \}, 
\quad
 j,  k \in {\Bbb Z}_+,
\end{align*}
one has 
\begin{align*}
(\gamma_j(H)\gamma_{n-j}&(H^k))_{u,w}
= \card \{ c_1\cdots c_k : c_l \in C, 1 \le l \le k,
s(c_1) = u, r(c_k)= w, \\
 & r(c_l) = s(c_{l+1}), 1 \le l \le k, 
 \ell(c_1) =j, \
\ell(c_1\cdots c_k) = n-j \}, 
\
1 \le j \le n, \ k \in {\Bbb Z}_+.
\end{align*}
It follows that
\begin{align*}
\log \zeta_{X_\C}(z) 
= &\sum_{n \in {\Bbb N}} \frac{z^n}{n} 
j \ \card\{  c_1\cdots c_k :c_l \in C, 1 \le l \le k,
r(c_l) = s(c_{l+1}),  1 \le l \le k,\\ 
 &  s(c_1) = r(c_k), \ 
\ell(c_1) =j, \
\ell(c_1\cdots c_k) = n  \} \\
= & \sum_{n \in {\Bbb N}} \frac{z^n}{n}
\sum_{1 \le j \le n}j \ \trace (\sum_{k \in {\Bbb Z}_{+}}\gamma_j(H)
\gamma_{n-j}(H^k) )\\
= &\sum_{n \in {\Bbb N}} \frac{z^n}{n}
\sum_{1 \le j \le n}j \ \trace (\gamma_j(H) \gamma_{n-j}((I-H)^{-1}))
 \\
= &\sum_{n \in {\Bbb N}} \frac{z^n}{n}
\sum_{0 \le j < n}(j+1) \trace (\gamma_{j+1}(H)
\gamma_{n-1-j}((I-H)^{-1})) \\
= & \sum_{n \in {\Bbb N}} \frac{z^n}{n}
\sum_{0 \le j < n} \trace (\gamma_{j}(H') \gamma_{n-1-j}(I-H)^{-1})) \\
= &\sum_{n \in {\Bbb N}} 
\frac{1}{n} \ \trace( \gamma_{n-1}(H' (I-H)^{-1}) z^n )\\
= &-\sum_{n \in {\Bbb N}} 
\trace( \gamma_{n}(\log (I-H)) z^n )\\
= & -\trace (\log (I-H)). 
\end{align*}
It is (see \cite[Section 1.1.10]{GJ}),
$$
\trace(\log ( I-H)) =  \log \det(I-H),
$$
and the theorem follows.
\end{pf}

We state Keller's formula for the case of 
a circular code (see \cite{Sta}, \cite{P}, 
and references given in \cite{BBG}) as a
corollary.

\begin{cor}\label{cor: 2.2}
For a circular code $\C$ 
$$
\zeta_{X_\C}(z) = \frac{1}{1-g_\C (z)}.
$$
\end{cor}
Note that also the formula for the zeta function of a subshift of finite type
in terms of a presenting polynomial matrix \cite{B} is also a special case 
of Keller's formula.

Let $G$ be a finite directed graph with adjacency matrix $A_G$.
We introduce the Markov-Dyck codes 
$\C_v, v\in {\cal V}$
of words $c = (c_k)_{1\le k \le K}$
$
\prod_{1 \le k \le K} c_k =P_v, \
\prod_{1 \le j \le J} c_k \ne P_v, 1 \le J <K. 
$
Standard methods of combinatorics (as for instance described in \cite{D})
give
$$
g_{\C_u}(z) = z^2 \sum_{v \in V} \frac{A_G(u,v)}{1 - g_{\C_v}(z)}, 
\qquad u \in 
{\cal V} 
$$
and by the implicit function theorem   
(1) has a unique solution
(see for instance \cite{Ma3},\cite{Ku}).
Set 
$$
\C = \cup_{v \in {\cal V}}\C_v.
$$
Also denote by 
$\C^-$ the set of admissible words that are concatenations of
an element 
(possibly empty)
of ${\cal P}^-$  
with a word in $\C$
and denote by 
$\C^+$ the set of admissible words that are concatenations of 
a  word in  $\C$
and an element  (possibly empty) of ${\cal P}^+$.
$(\C^-,s, r)$ and $(\C^+ , s, r) $ are circular Markov codes.
Denote by 
$D(A_G,z) $ the diagonal matrix with entries $g_{\C_v}(z), v \in \V$,
and denote by $D^*(A_G,z)$ the diagonal matrix 
with entries
$g_{\C^*_v}(z) =\frac{1}{1-g_{\C_v}(z)}, v \in {\cal V}$.
\begin{thm}\label{thm:2.3}
The zeta function of the Markov Dyck shift $D_G$ is
\begin{align*}
\zeta_{D_G}(z) 
& = \frac{1}{\det(( I - D(A_G,z) - A_G z)( I - D^*(A_G,z) A_G z))}\\
& = \frac{\det(D^*(A_G,z))}{\det(( I - D^*(A_G,z) A_G z))^2}
\end{align*}
\end{thm}
\begin{pf}
Since
$ 
\sum_{ k \in {\Bbb Z}_+ } A_G^k z^k = (I - A_G z)^{-1} 
$
one has
$$
H^{(\C^+)}(z) = D(A_G,z)(I - A_G z)^{-1},
$$
and $H^{(\C^-)}(z)$ is the adjoint of $H^{(\C^+)}(z)$.
Applying Proposition 2.1 and collecting
 all contributions to the zeta function,
one has
\begin{align*}
\zeta_{D_G}(z)
& = (\prod_{u \in \V}g_{\C^*_u}(z)^{-1}) 
   \det(I - A_G z)^{-2}\det(I - D(A_G,z) (I - A_G z)^{-1})^{-2}\\
& = (\prod_{u \in \V}g_{\C^*_u}(z)^{-1})\det({D^*(A_G,z)}^{-1} - A_G z )^{-2}\\
& = \det(({D^*(A_G,z)}^{-1} - A_G z)( I - D^*(A_G,z) A_G z))^{-1}\\
& = \frac{1}{\det(( I - D(A_G,z) - A_G z)( I - D^*(A_G,z) A_G z))}\\
& = \frac{\det({D^*(A_G,z)})}{\det( I - D^*(A_G,z) A_G z )^2}.
\end{align*}
\end{pf}
Inserting into the formula for the case of 
the graph with one vertex and $N$-loops,
the generating function 
$$
g_{\C_v}(z) = \frac{ 1 - \sqrt{1-4Nz^2}}{2}
$$
one obtains again the zeta function of the Dyck shift $D_N$ as
$$
\zeta_{D_N}(z) = \frac{2(1 + \sqrt{1-4Nz^2})}{(1-2Nz + \sqrt{1-4Nz^2})^2}
$$
( see \cite{Ke}).

\section{Topological entropy}

\begin{prop}\label{prop:3.1}
For the Markov-Dyck shift $D_G$
$$
h(D_G) = \lim_{n \to \infty} \frac{1}{n}\log \Pi_n(D_G).
$$
\end{prop}
\begin{pf}
A word $b \in {\cal L}  _n(D_G), n \in {\Bbb N},$
determines words
$a^+(b) \in {\cal D}^+, a^-(b) \in {\cal D}^-,$
by
$$
\prod_{1 \le m \le n}b_m = a^+(b) a^-(b),
$$
as well as indices
$I^-(b), I^+(b)$, $1 \le I^-(b), I^+(b) \le n,
$
by
$$
I^-(b)  = \min \{ i : \prod_{1 \le j \le i}b_j = a^+(b) \},\qquad
I^+(b)  = \max \{ i : \prod_{i \le j \le n}b_j = a^-(b) \}.
$$
Denote by
${\cal K}_n$ 
the set of $b \in {\cal L}(D_G)$ of length $n \in {\Bbb N}$ 
such that $I^-(b) =1$.
Choose for $u,w \in \V$
a path $c(u,w)$ in $G$ from $u$ to $w$ 
of shortest length $\lambda(u,w)$ and set 
$$
L = \max_{u,w \in \V} \lambda(u,w).
$$
We define a mapping
$\Psi_n$ of ${\cal L}_n(D_G)$ into
$\cup_{n \le m \le n+2L}{\cal K}_m$
by
$$
\Psi_n(b) = b_{[I^{-}(b),n]})*c(b_n,b_{I^{-}(b)})
*a^-*c(b_1,b_{I^{-}(b)}),
\qquad
a^+(b) = a^+, b \in {\cal L}_n(D_G), n \in {\Bbb N}.
$$
A word in ${\cal K}_m$, $n \le m \le m + 2L,$
has at most $(n+2L)\card {\cal V}$ inverse images under the mapping 
$\Psi_n$, therefore
\begin{equation}
\card \ {\cal L}_n(D_G)  \le ( n+2L) \sum_{n \le m \le n+2L} |{\cal K}_m|. 
\label{eqn:3.1}
\end{equation}
Every word in ${\cal K}_m, n \le m \le n+2L,$
determines a periodic point in $D_G$
and \eqref{eqn:3.1} implies that 
$$
\lim_{n \to \infty}\frac{1}{n} \log |{\cal L}_n(D_G)|
\le 
\liminf_{n \to \infty} \frac{1}{n} \log \Pi_n(D_G).
$$
\end{pf}
\begin{cor}\label{cor:3.2}
For the Markov-Dyck shift $D_G$,
the topological entropy 
$
h(D_G)
$
is the minimum positive solution of the equation:
$$
\det( I - D^*(A_G,z) A_G z) =0.
$$
\end{cor}

For $v \in {\cal V}$ let $X_v$ denote the subsystem of the Markov-Dyck shift 
$D_G$
that is obtained by excluding the words 
$$
e(-)e(+),\qquad  e \in E, \quad r(e) \in \V \backslash\{v\}
$$
and the words
$$\
f(+)g(-),\qquad  g,f \in E, \quad s(f) \in \V\backslash\{v\}.
$$
We will estimate the asymptotic grooth rate
of the periodic points of $X_v$ which, by a proof that is similar to the
proof of Proposition \ref{prop:3.1}, 
is actually equal to the topological entropy of $X_v$.
In this way, we will also obtain estimates of the topological entropy of the
Markov-Dyck shifts.

For $v \in {\cal V}$ denote by ${\cal D}_v$ the circular code of elementary
Markov-Dyck words
that start and end at $v$,
and denote by $\C_v$ the circular code of paths in $G$
that start at $v$ and end at $v$ when they return for the first time to $v$.
$\rho$ denotes the inverse of the Perron eigenvalue of $A_G$.
We denote by $p(z)$ the determinant of the matrix 
$ I - A_G z$ 
and
by $p_v(z)$
the determinant of the matrix
$ I - A_G z$ 
with the $v$-th row and the $v$-th column deleted, 
$v \in {\cal V}$.
We set 
$$
q_v = \frac{p}{p_v},\qquad v \in {\cal V}.
$$
\begin{prop}\label{prop:3.3}
$$
g_{{\cal D}_v^* * \C_v^*}(z) = \frac{ 1 - \sqrt{1 - 4
g_{\C_v}(z^2)}}{2(1-g_{\C_v}(z))},
\qquad
v \in {\cal V}.
$$
\end{prop}
\begin{pf} 
One has 
$g_{{\cal D}_v * \C_v^*} = g_{{\cal D}_v} g_{\C_v^*}$
and
$g_{{\cal D}_v}$ satisfies the equation
$$
g_{{\cal D}_v}(z) = \frac{g_{\C_v}(z^2)}{ 1-g_{{\cal D}_v}(z)}.
$$
It follows that
$$
g_{{\cal D}_v}(z) = \frac{1}{2} (1 - \sqrt{1 - 4 g_{\C_v}(z^2)}),
$$
which yields the proposition.
\end{pf}
\begin{thm}\label{thm:3.4}
Let $v \in {\cal V}$ be such that 
\begin{equation} 
q_v(\rho^2) > \frac{3}{4}. 
\label{eqn:3.2}
\end{equation}
Then 
\begin{equation}
h(X_v) >
- \log \rho +
\frac{
p_v(\rho)[q_v(\rho^2)-\frac{3}{4}-\frac{1}{2}\sqrt{
q_v(\rho^2)-\frac{3}{4}}]
}{
\rho p'(\rho)[\rho + \sqrt{q_v(\rho^2)-\frac{3}{4}}]
}.
 \label{eqn:3.3}
\end{equation}
\end{thm}
\begin{pf}
Corollary 2.2 and Proposition 3.3
 imply that $h(X_v)$ is 
equal to $-\log \kappa$, where
$ \kappa$ is the solution of the equation 
$$
1 = 2 g_{\C_v}(z) - \sqrt{1 - 4 g_{\C_v}(z^2)}, \qquad 0 < z < \rho,
$$
or, equivalently,
of the equation
$$
\frac{1}{2} = q_v(z) + \sqrt{q_v(z^2)-\frac{3}{4}}, \qquad 0 < z < \rho.
$$
One has the estimate
$$
\rho -\kappa > \frac{
p_v(\rho)[q_v(\rho^2)-\frac{3}{4}-\frac{1}{2}\sqrt{
q_v(\rho^2)-\frac{3}{4}} }{p'(\rho)[\rho + \sqrt{q_v(\rho^2)-\frac{3}{4}}]}
$$
and \eqref{eqn:3.3} follows.
\end{pf}
\begin{thm}\label{thm:3.5}
Let $v \in {\cal V}$ be such that 
\begin{equation} 
q_v(\rho^2) < \frac{3}{4}. 
\label{eqn:3.4}
\end{equation}
Then 
\begin{equation}
h(X_v) >
-\log \rho + 
\frac{ p_v(\rho)}{
2\rho^2p'(\rho)
}[q_v(\rho^2)-\frac{3}{4}]
 \label{eqn:3.5}
\end{equation}
\end{thm}
\begin{pf}
Corollary 2.2 and Proposition 3.3 
and \eqref{eqn:3.4} imply that 
$h(X_v)$ is greater than or equal to $-\log \kappa$,
where
$\kappa$ is the solution of the equation
$$
1- q_v(z^2) = \frac{1}{4}, \qquad 0 < z < \rho.
$$
One has the estimate
$$
\rho - \kappa >
\frac{p_v(\rho)}{2 p'(\rho)}[q_v(\rho^2) - \frac{3}{4}] 
$$ 
and \eqref{eqn:3.4} follows.
\end{pf}
\begin{lem}\label{lem:3.6}
Let 
\begin{equation}
\rho < \frac{1}{4}.
\label{eqn:3.6}
\end{equation}
Then there is a $v \in {\cal V}$ such that
$$
q_v(\rho^2) > \frac{3}{4}.
$$
\end{lem}
\begin{pf}
Assume the contrary.
Then 
$$
\frac{3}{4} \card {\cal V}  
\le 
\sum_{v \in \V} \frac{p_v(\rho^2)}{p(\rho^2)}
= \trace ({I-A_{\rho^2}}^{-1})
 \le \frac{1}{1-\rho} 
\card {\cal V},
$$
contradicting \eqref{eqn:3.6}.
\end{pf}

Theorems \ref{thm:3.4} and \ref{thm:3.5} 
have as a corollary an estimate for of the topological
entropy of the Markov-Dyck shifts. We state here a corollary of 
Theorem \ref{thm:3.4}.
\begin{cor}\label{cor:3.7}
Let 
$$
\rho < \frac{1}{4}.
$$
Then
$$
h(D_G)  > - \log \rho       + \max_{ \{ v \in {\cal V}:q_v(\rho^2) > \frac{3}{4}
\}} \frac{
p_v(\rho)[q_v(\rho^2)-\frac{3}{4}-\frac{1}{2}\sqrt{
q_v(\rho^2)-\frac{3}{4}}]
}{
\rho p'(\rho)[\rho + \sqrt{q_v(\rho^2)-\frac{3}{4}}]
}.
$$
\end{cor}
\begin{pf}
The corollary follows from Theorem \ref{thm:3.4} 
by means of Lemma \ref{lem:3.6}.
\end{pf}

\section{A class of examples}
We consider first the Fibonacci-Dyck shift $D_F$ 
that is produced by the directed graph (Figure 1)
with adjacency matrix 
$
F =
\begin{bmatrix}
1 & 1 \\
1 & 0 
\end{bmatrix}.
$
Here
\begin{figure}[htbp]
\begin{center}
\unitlength 0.1in
\begin{picture}( 21.9100,  6.9000)( 16.5000,-24.0600)
%
\special{pn 8}%
\special{ar 2092 2064 78 66  0.0000000 6.2831853}%
%
\special{pn 8}%
\special{ar 3762 2058 80 66  0.0000000 6.2831853}%
%
\special{pn 8}%
\special{ar 1832 2058 182 154  0.3190696 6.0756891}%
%
\special{pn 8}%
\special{pa 2002 2112}%
\special{pa 2006 2106}%
\special{fp}%
\special{sh 1}%
\special{pa 2006 2106}%
\special{pa 1958 2158}%
\special{pa 1982 2154}%
\special{pa 1994 2176}%
\special{pa 2006 2106}%
\special{fp}%
%
\special{pn 8}%
\special{ar 2946 2128 818 278  6.2480092 6.2831853}%
\special{ar 2946 2128 818 278  0.0000000 3.1562795}%
%
\special{pn 8}%
\special{pa 3762 2122}%
\special{pa 3762 2118}%
\special{fp}%
\special{sh 1}%
\special{pa 3762 2118}%
\special{pa 3742 2186}%
\special{pa 3762 2172}%
\special{pa 3782 2186}%
\special{pa 3762 2118}%
\special{fp}%
%
\special{pn 8}%
\special{pa 2134 2006}%
\special{pa 2136 1974}%
\special{pa 2148 1946}%
\special{pa 2166 1918}%
\special{pa 2188 1896}%
\special{pa 2214 1876}%
\special{pa 2240 1858}%
\special{pa 2268 1842}%
\special{pa 2296 1830}%
\special{pa 2326 1816}%
\special{pa 2356 1804}%
\special{pa 2386 1794}%
\special{pa 2416 1784}%
\special{pa 2448 1776}%
\special{pa 2478 1768}%
\special{pa 2510 1760}%
\special{pa 2542 1754}%
\special{pa 2572 1748}%
\special{pa 2604 1742}%
\special{pa 2636 1738}%
\special{pa 2668 1734}%
\special{pa 2700 1730}%
\special{pa 2732 1728}%
\special{pa 2762 1724}%
\special{pa 2794 1722}%
\special{pa 2826 1720}%
\special{pa 2858 1720}%
\special{pa 2890 1718}%
\special{pa 2922 1716}%
\special{pa 2954 1716}%
\special{pa 2986 1718}%
\special{pa 3018 1718}%
\special{pa 3050 1718}%
\special{pa 3082 1720}%
\special{pa 3114 1722}%
\special{pa 3146 1726}%
\special{pa 3178 1728}%
\special{pa 3210 1730}%
\special{pa 3242 1736}%
\special{pa 3274 1740}%
\special{pa 3306 1744}%
\special{pa 3336 1750}%
\special{pa 3368 1756}%
\special{pa 3400 1764}%
\special{pa 3430 1770}%
\special{pa 3462 1778}%
\special{pa 3492 1788}%
\special{pa 3522 1796}%
\special{pa 3554 1806}%
\special{pa 3582 1820}%
\special{pa 3612 1832}%
\special{pa 3640 1846}%
\special{pa 3668 1862}%
\special{pa 3694 1880}%
\special{pa 3718 1902}%
\special{pa 3740 1926}%
\special{pa 3758 1952}%
\special{pa 3768 1982}%
\special{pa 3768 2000}%
\special{sp}%
%
\special{pn 8}%
\special{pa 2134 1992}%
\special{pa 2134 2006}%
\special{fp}%
\special{sh 1}%
\special{pa 2134 2006}%
\special{pa 2154 1940}%
\special{pa 2134 1954}%
\special{pa 2114 1940}%
\special{pa 2134 2006}%
\special{fp}%
\put(20.9000,-20.7000){\makebox(0,0){1}}%
\put(37.6000,-20.5000){\makebox(0,0){2}}%
\end{picture}%
\end{center}
\caption{}
\end{figure}

\begin{align}
g_{\C_1}(z) & = (g_{\C_1^*}(z) + g_{\C_2^*}(z))z^2,  \\
g_{\C_2}(z) & =  g_{\C_1^*}(z)z^2, 
\end{align}
where 
\begin{equation}
g_{\C_1^*}  = \frac{1}{ 1 - g_{\C_1}}, \qquad 
g_{\C_2^*}  = \frac{1}{ 1 - g_{\C_2}}, 
\end{equation}
or
\begin{equation}
g_{\C_1}  = 1 - \frac{1}{g_{\C_1^*}}, \qquad 
g_{\C_2}  = 1 - \frac{1}{ g_{\C_2^*}}.
\end{equation}
From (4.1), (4.3) and (4.4)
\begin{equation}
g_{\C_1^*}(z) = 1 +  g_{\C_1^*}(z)(g_{\C_1^*}(z) + g_{\C_2^*}(z)) z^2, 
\end{equation}
and from (4.2) and (4.3)
\begin{equation}
g_{\C_2^*}(z) = 1 + g_{\C_1^*}(z)g_{\C_2^*}(z)z^2. 
\end{equation}
From (4.2), (4.4) and (4.5)
\begin{equation}
 g_{\C_2^*}(z)^3 z^2- g_{\C_2^*}(z) + 1 = 0, 
\end{equation}
and from (4.6) and (4.7) 
\begin{equation}
g_{\C_1^*} =  g_{\C_2^*}^2. 
\end{equation}
From  (4.2) and (4.3)
\begin{equation}
 \det{(I - Fz - D(F,z))} 
  = \frac{z}{g_{\C_2}}(
g_{\C_2}(z)^2 - (2 z +1)g_{\C_2}(z) + z ).
\end{equation}
Setting 
$\xi(z) =  g_{\C_2^*}(z) z$  one has from (4.7)
\begin{equation}
\xi(z)^3 - \xi(z) + z =0    
\end{equation}
and from (4.9) and (4.10)
\begin{equation}
 \det{(I - Fz - D(F,z))} 
  = - \frac{z^2}{\xi(z)^2}(
2 \xi(z)^2 + \xi(z) -1 ).
\end{equation}
By Theorem 2.3
and by (4.4) and (4.11)
\begin{equation}
\zeta_{D_F}(z) = \frac{\xi(z)}{z( 2\xi(z)^2 + \xi(z) - 1)^2},
\end{equation}
where one identifies $\xi$ as the solution 
of equation (4.10) vanishing at the origin that is given by 
\begin{equation}
\xi(z) = \frac{2}{\sqrt{3}} 
\sin(\frac{1}{3}\arcsin\frac{3\sqrt{3}}{2} z), 
\qquad 0 \le z \le
\frac{2}{3 \sqrt{3}}. 
\end{equation}
By Theorem 3.1
and by (4.12)
the topological entropy of the Fibonacci-Dyck shift is equal to the negative
logarithm of the solution of
$$
2 \xi(z)^2 + \xi (z) - 1 = 0.
$$
By (4.10) (or by (4.13)),
\begin{equation}
h(D_F) =   3 \log 2 - \log 3.
\end{equation}

\medskip


We turn to the Markov-Dyck shift 
that is produced by the directed graph with adjacency matrix
$
F(a,b,c) = 
\begin{bmatrix}
a & b \\
c & 0 
\end{bmatrix},
a,b,c \in {\Bbb N}.
$
Here
\begin{align}
g_{\C_1}(z) & = (a g_{\C_1^*}(z) + b  g_{\C_2^*}(z)) z^2, \\
g_{\C_2}(z) & = c g_{\C_1^*}(z) z^2 
\end{align}
and one has from (4.16) that 
\begin{equation}
g_{\C_1}(z)  = 1 - \frac{ c z^2}{ g_{\C_2}(z)}. 
\end{equation}
From (4.15) and (4.17)
\begin{equation}
a g_{\C_2}(z)^3 - ( a + c) g_{\C_2}(z)^2 +
c( 1 +(c - b)z^2) g_{\C_2}(z) -c^2 z^2 =0.
\end{equation}
From (4.17)
\begin{equation}
 \det{(I - F(a,b,c)z - D(F(a,b,c),z))} 
  = \frac{z}{g_{\C_2}}(
a g_{\C_2}(z)^2 - (a + c(1 + b) z )g_{\C_2}(z) + c z ).
\end{equation}
Theorem 2.3
and  (4.17) and (4.19)
give
\begin{equation}
\zeta_{D_{F(a,b,c)}}(z) 
= \frac{c g_{\C_2}(z)(1 - g_{\C_2}(z))}{(
a g_{\C_2}(z)^2 - ( a + c(1 + b)z) g_{\C_2}(z) + c z)^2}.
\end{equation}
Setting 
\begin{align*}
\mu(z) &=  ( c -a )^2 + ac - 3ac ( c -b ) z^2, \\
\nu(z) &=  2( a + c)^3 - 9ac( a+b) - ( c - b + 27 a^2 c^2) z^2, 
\end{align*}
one identifies  
$g_{\C_2}(z)$
as the solution of (4.18)
that vanishes at the origin
\begin{equation*}
g_{\C_2}(z) = 
- \frac{a+c}{3a} +
\frac{2}{3 a} \sqrt{\mu(z)} 
\cos(\frac{1}{3}( 2 \pi + \arccos\frac{\nu(z)}{\mu(z) \sqrt{\mu(z)}})) 
\end{equation*}
(For the case $ a = b = c= 1$ compare here (4.13)).

We determine the topological entropy of 
$D_{F(a,b,c)}, \ a, b, c, \in {\Bbb N}$.
Set
\begin{align*}
P_{a,b,c}(z)  
=  & (1 +c) [ a(b-c) - c (1 + b)^2] z^3 \\ 
+  & ( c [(1 + b) ( 1 + c) -2ab] + a( 1 + a -b) ) z^2 \\ 
+ & ( bc - a - (1+a)(a-c) ) z 
+ a - c.
\end{align*}

\begin{thm}
\qquad
\qquad
\qquad
\qquad
 \begin{enumerate}
\renewcommand{\labelenumi}{(\alph{enumi})}
\item  $h(D_{F(a,b,c)})$
is equal to the negative logarithm of the smallest positive solution 
of $P_{a,b,c}(z) =0, \ a, b ,c \in {\Bbb N}.$
\item  $h(D_{F(a,b,a+b)}) = - \log ( 1 + a + b), \ a, b \in {\Bbb N}.$
\end{enumerate}
\end{thm}
\begin{pf}
Let $ z >0$ be such that the equations 
\begin{align}
a y^2 & - ( a + c(1 + b) z) y + c z = 0 \\
\intertext{and}
a y^3 & - ( a + c) y^2 + c(1 + ( c - b) z^2) y  - c^2 z^2 = 0 
\end{align}
have a common solution $y$.
Then $y$ also solves the equation 
\begin{equation}
( 1 - ( 1 + b)z) y^2 - ( 1 - z + ( c -b) z^2) y + c z^2 = 0
\end{equation}
and, as is seen from (4.21) and (4.23),
it also solves the equation
\begin{equation}
( 1 - ( 1 +a +  b)z) y= 1 - ( 1 + a) z - b( 1 + c) z^2.
\end{equation}
From (4.21) and (4.24)
\begin{equation}
\begin{split}
 b z P_{a,b,c}(z) 
 =& ( 1 - ( 1 + a) z - b ( 1 + c)z^2) \\
  & \{ a( 1 - ( 1 + a) z - b ( 1 + c)z^2)
  - ( 1 - ( 1 + a+ b) z) ( a + c(1 +b)z) \}\\ 
  & + c z ( 1 - ( 1 + a + b)z)^2 = 0.  
\end{split}
\end{equation}
This shows that for every $ z >0$ such that equations (4.21) and (4.22) 
have a common solution, $P_{a,b,c}(z) =0$.

Equation (4.23) is a multiple of equation (4.21) precisely if
$ c = a +b$ and
$ z = \frac{1}{1 + a+b}$,
and from this one sees, consulting (4.24),
that both solutions of equation (4.21) 
are also solutions of equation (4.22)
precisely if 
$ c = a +b$
and
$ z = \frac{1}{1 + a+b}$.
Moreover,
as is seen from (4.25),
$P_{a,b,c}(z) =0$ has the root 
$ \frac{1}{1 + a+b}$
if and only if
$ c = a +b$.

For the case that $ c \ne a+b$, let $z >0$,
\begin{equation}
P_{a,b,c}(z) =0,
\end{equation}
and reverse the argument, setting
\begin{equation}
y = \frac{ 1 - ( 1 + a)z - b(1 + c) z^2}{1 - ( 1 + a + b) z}. 
\end{equation}
Consult  then (4.24) and find from (4.26)
that $y$ as given by (4.27)  solves equation (4.21) 
and therefore also 
equations (4.23) and (4.22).
Apply now Theorem 2.3 in conjuction with Theorem 2.1 
together with (4.18) and (4.20)
to prove part (a) of the theorem for the case $c \ne a+b$.

Consider the case that $c = a+b$.
One checks that $ \frac{1}{1 + a + b}$ is the unique positive root of
$P_{a,b,a+b}(z) = 0$.
It has already been shown that 
$P_{a,b,a+b}(z) = 0$ 
for every $ z>0$ 
such  that equations (4.21) and (4.22), or, in this case, 
the equations
\begin{equation}
a y^2  - ( a + ( a +b)(1 + b) z) y +(a+b) z = 0 
\end{equation}
and
\begin{equation*}
a y^3  - ( 2a + b) y^2 + (a+b)(1 + a z^2) y  - (a+b)^2 z^2 = 0 
\end{equation*}
have a common solution.
One checks that
for
$ z = \frac{1}{ 1+ a+b}$ 
a root of (4.28),
in fact
the smaller one,
 is equal to 
$g_{\C_2}(\frac{1}{1+ a+ b})$.
Apply now again  Theorem 2.3 in conjuction with 
Theorem 2.1 together with
(4.18)
and (4.20)
to prove part (a) of the theorem for the case that 
$c= a+b$ and also part (b).
\end{pf}
The corollary reconfirms (4.14).
\begin{cor}
\begin{equation*}
h(D_{F(a,1,a)}) =
\log(a+1) - \log(a +2) + \log(a+3), \ a \in {\Bbb N}
\end{equation*}
\end{cor}
Determining the topological entropy of 
$D_{F(a,b,a)}, \ a,b \in {\Bbb N}$
and of 
$D_{F(a,b,c)}$ where
$a,b,c \in {\Bbb N}$
solve
$a(b-c) - c(1+b)^2 =0$,
reduces to solving a quadratic equation.
In the various further cases
expressions for the topological entropy are obtained from formulas
of
Tartaglia and Vieta.


\begin{thebibliography}{99}


\bibitem{AA}
{\sc G.~Abrams and G.~Aranda-Pino},       
{\it The Leavitt path algebra of a graph},
J. Algebra \
{\bf 293}(2005), pp.\ 319--334.



\bibitem{BP}{\sc J.~Berstel and D. ~Perrin},
{\it Theory of codes}, Academic Press, London
(1985).



\bibitem{B}
{\sc  M.~Boyle},
{\it Symbolic dynamics and matrices},
 Combinatorial and Graph-theoretical Problems in Linear Algebra, 
 IMA Volumes in Mathematics and its Applications 
{\bf 50}(1993), pp.\ 1--38.


\bibitem{BBG}
{\sc M.~Boyle, J.~Buzzi and R.~Gomez},       
{\it Almost isomorphism for countable state Markov shifts},
J. reine angew. Math. \
{\bf 592}(2006), pp. \ 23--47.






\bibitem{C}{\sc J. ~Cuntz}, 
{\it Simple $C^*$-algebras generated by isometries},
Commun.\ Math.\ Phys.\
{\bf 57}(1977), pp.\ 173--185.



\bibitem{CK}{\sc J. ~Cuntz and W. ~Krieger},
{\it A class of $C^*$-algebras and topological Markov chains},
 Inventions Math.\
 {\bf 56}(1980), pp.\ 251--268.


\bibitem{D}
{\sc M.~Delest},
{\it Algebraic languages: a bridge between combinatorics and 
computer science},
 Dimacs \
{\bf 24}(1996), pp.\ 71--87.






\bibitem{EW}
{\sc M.~Enomoto and Y.~Watatani},
{\it  A graph theory for C*-algebras},
Math. Japon.\
{\bf 25}(1980), pp.\ 435-- 442.

\bibitem{GJ}
{\sc I.~P.~Goulden and D.~M.~Jackson},
{\it Combinatorial Enumeration},
John Wiley, \
New York, 1983.


 
\bibitem{HIK}{\sc T. ~Hamachi, K. ~Inoue and W. ~Krieger}, 
{\it  Subsystems of finite type and semigroup invariants of subshifts},
 preprint.                      



\bibitem{In}{\sc  K. ~Inoue }, 
{\it  The zeta function, periodic points and entropies of the Motzkin shift},
preprint, arXiv:math.DS/0602100.   

\bibitem{Ke}{\sc G. ~Keller},
{\it  Circular codes, loop counting, and zeta-functions}, 
J.\ Combinatorial Theory
{\bf  56}
(1991),
pp.\ 75--83.




\bibitem{Ki}{\sc B.~P. ~Kitchens},
{\it Symbolic dynamics}, Springer-Verlag, Berlin, Heidelberg and New York
(1998).



\bibitem{Kr1}{\sc W. ~Krieger},
{\it  On the uniqueness of the equilibrium state},  Math.\ Systems Theory
{\bf  8}
(1974),
pp.\ 97--104.



\bibitem{Kr2}
{\sc W. ~Krieger},
{\it  On a syntactically defined invariant of symbolic dynamics},
 Ergodic Theory Dynam. Systems
{\bf 20}(2000), pp.\ 501--516.


\bibitem{Kr3}{\sc W. ~Krieger},
{\it  On subshifts and semigroups},  Bull.\ London Math.
{\bf  38}
(2006),
pp.\ 617--624.




 
\bibitem{KM}
{\sc W. ~Krieger and K. ~Matsumoto},
{\it A lambda-graph system for the Dyck shift and its K-groups},                  Doc. Math.\
{\bf 8}(2003), pp. \ 79--96.




\bibitem{Ku}
{\sc W. ~Kuich},
{\it On the entropy of context-free languages},
Information and Control
{\bf 16}(1970),pp.\ 173--200.

\bibitem{L}
{\sc W.~G.~Leavitt},
{\it The module type of homomorphic images},
Duke Math.J.
{\bf 32}(1965), pp.\  305--31.




\bibitem{LM}{\sc D. ~Lind and B. ~Marcus},
{\it An introduction to symbolic dynamics and coding},
 Cambridge University Press, Cambridge
(1995).







\bibitem{Ma1}{\sc K. Matsumoto},
{\it  A simple purely infinite 
$C^*$-algebra associated with a lambda-graph system of Motzkin shift},  Math.\ Z.\ {\bf  248}(2004), pp.\ 369--394.




\bibitem{Ma2}{\sc K. Matsumoto},
{\it  K-theoretic invariants and conformal measures of the Dyck subshift}, 
Internat.\ J.\ Math.\ {\bf 16} (2005), pp.\ 213--248.


%


\bibitem{Ma3}{\sc K. Matsumoto},
{\it  Cuntz-Krieger algebras and a generalization of Catalan numbers},
preprint, arXiv:math.OA/0607517.


\bibitem{Ma4}{\sc K. Matsumoto},
{\it  
$C^*$-algebras arising from Dyck systems of topological Markov chains}, 
preprint, arXiv:math.OA/0607518.

\bibitem{Ma5}{\sc K. Matsumoto},
{\it  K-theory for the simple $C^*$-algebra of the Fibonacci Dyck system},
preprint, arXiv:math.OA/0607519.

\bibitem{NP}{\sc M. ~Nivat and J.-F. ~Perrot},
{\it  Une g\'{e}n\'{e}ralisation du mono\^{i}de bicyclique},
 C. R. Acad. Sc. Paris, 
 {\bf 271}
(1970), pp.\ 824--827



\bibitem{Pa}
{\sc L.~T.~Paterson}
{\it Graph inverse semigroups, groupoids and their C*-algebras},
J. Operator Theory
{\bf 48}(2004), pp. \ 645--662.

\bibitem{P}
{\sc D. ~Perrin},
{\it Algebraic combinatorics on words},
Algebraic Combinatorics and Computer Science, H.Crapo and 
G.-C.Rota, Eds. 
 Springer  2001, pp. \ 391--430.




\bibitem{Sta}{\sc R. P. ~Stanley},
{\it Enumerative combinatrics  I},
Wadsworth $\And$ Brooks/Cole Advanced Books $\And$ Software, Monterey,CA,
(1986).




\end{thebibliography}
\end{document}